\documentclass[12pt,a4paper]{article}

\usepackage{graphicx}
\usepackage{amsthm}
\usepackage{amsmath}
\usepackage{amssymb}

\newtheorem{remark}{Remark}[section]

\newtheorem{theorem}[remark]{Theorem}
\newtheorem{propos}[remark]{Proposition}
\newtheorem{corol}[remark]{Corollary}
\newtheorem{lemma}[remark]{Lemma}
\newtheorem{defin}[remark]{Definition}
\theoremstyle{definition}
\newtheorem{pict}[remark]{Figure}

\newcommand{\R}{R}
\newcommand{\Rn}{{\R}^n}
\newcommand{\Comp}{K(\Rn)}

\newcommand{\inK}{\in\Comp}

\newcommand{\PrjXonY}[2]{\Pi_{#2}{(#1)} }
\newcommand{\Pair}[2]{\Pi({#1},{#2})}
\newcommand{\Bspline}[2]{b_#1\left(#2\right)}

\newcommand{\Bern}[1]{B_N(#1, x)}
\newcommand{\MA}[3]{#1 \oplus_{\,#3} #2}
\newcommand{\Sum}[2]{\sum_{#1}^{#2}}
\newcommand{\Map}[1]{F:#1 \rightarrow \Comp }

\newcommand{\Lip}[1]{Lip\,(#1, {\cal L}) }

\begin{document}
\title {Approximations of Set-Valued Functions by Metric Linear Operators}
\bigskip

\author{Nira Dyn, Elza Farkhi, Alona Mokhov \\
    School of Mathematical Sciences \\
    Tel-Aviv University, Israel}
\bigskip

\date{}
\maketitle \hyphenation{set-valued single-valued ope-ra-tor
ope-ra-tors multi-function multi-functions pa-pa-me-tri-za-tion }

{ \small {\normalsize \textbf {Abstract.}} In this work, we
introduce new approximation operators for univariate set-valued
functions with general compact images. We adapt linear
approximation methods for real-valued functions by replacing
linear combinations of numbers with new metric linear combinations
of finite sequences of compact sets, thus obtaining "metric
analogues" operators for set-valued functions. The new metric
linear combination extends the binary metric average of Artstein.
Approximation estimates for the metric analogue operators are
derived. As examples we study metric Bernstein operators, metric
Shoenberg operators and metric polynomial interpolants. }
\medskip

\noindent{ \small {\normalsize \textbf{Key words:}} compact sets,
Minkowski linear combination, metric average, set-valued
functions, piecewise linear set-valued functions, selections,
linear approximation operators, Bernstein polynomial
approximation, Schoenberg spline approximation, polynomial
interpolation. }

\noindent{ \small {\normalsize \textbf{Mathematics Subject
Classification 2000:}} 26E25, 54C65, 41A35, 41A36}

\section {Introduction} \label{Sect_Intro}

In this work, we adapt approximating operators for real-valued
functions to set-valued functions (multifunctions, SVFs), by
replacing an operation between numbers with an operation between
sets. The known approximation methods, based on Minkowski sums of
sets, fail to approximate, when the images of a multifunction are
not convex. In case of Bernstein-type operators and subdivision
operators there is a phenomenon of "convexification"
(\cite{Vitale:Bern_Mink, Dyn-Farkhi:ConvRates}).

In~\cite{Artstein:MA} a binary operation between sets, the "metric
average", is introduced and the metric piecewise linear
interpolant based on it is shown to approximate continuous SVFs
with general images. The use of this operation in the adaptation
of known approximation methods to SVFs, requires a representation
of the approximation operators by repeated binary averages. Such a
representation exists for any operator which reproduces constants,
but is not unique~\cite{DynWal}. This non-uniqueness leads to
different operators and it is not clear what are the appropriate
adaptations. Spline subdivision schemes represented by repeated
averages \cite{Dyn-Farkhi:Subdiv_MA} and the Schoenberg operators
defined in terms of the de~Boor
algorithm~\cite{Dyn-Mokhov:Approx_MA} are proved to approximate
SVFs with general compact images. Yet, for the adaptation of the
Bernstein operators based on the de~Casteljau algorithm we could
obtain an approximation result only for SFVs with images in $\R$
all consisting of the same number of disjoint
intervals~\cite{Dyn-Mokhov:Approx_MA}.

In this paper we introduce a set-operation on a finite sequence of
compact sets, termed "metric linear combination", which extends
the metric average. We adapt approximation methods for real-valued
functions to SVFs, replacing linear combinations of numbers by the
metric linear combinations of sets. We prove that this adaptation
of any linear operator, approximating continuous real-valued
functions, approximates continuous SVFs of bounded variation. In
particular for Lipschitz continuous SVFs, sharper error estimates
are obtained. Approximation results for set-valued functions which
are only continuous, are given for a limited class of operators.
It should be noted that our adaptation method is not restricted to
positive operators. The approximation results are specialized to
the Shoenberg spline operators and the Bernstein polynomial
operators. Also the adaptation of polynomial interpolation to SVFs
is presented as examples of non-positive operators. This
adaptation is illustrated by two metric parabolic interpolants.

An outline of the paper is as follows. The next section contains
basic definitions, notation and known results. In
Section~\ref{Sect_Chains} we introduce the metric linear
combination between a finite number of ordered sets and define
metric linear operators for multifunctions based on it. In
Section~\ref{Sect_PLApprox} properties of the metric piecewise
linear intorpolant are considered. In particular a representation
of it by a specific set of selections is studied. Similar
selections are used in~\cite{HenryHermes:ContSelections}
and~\cite{Dommisch:LipSelections} to prove the existence of a
continuous selection for a continuous SVF of bounded variation,
and of a representation of a Lipschitz SVF respectively. In
Section~\ref{Sect_ApproxByOperators} we derive approximation
results for the metric linear approximation operators, based on
the results in Section~\ref{Sect_PLApprox}. Finally, in
Section~\ref{Sect_Examples}, we specialize these results to some
classical approximation operators.

\section {Preliminaries}
\label{Sect_Prelim}

First we present some definitions and notation.

\noindent$\bullet$ $\Comp$ is the collection of all compact
nonempty subsets of~$\Rn$.

\noindent$\bullet$ A linear Minkowski combination of two sets $A$
and $B$ from $\Comp$ is
  $$\lambda A + \mu B = \{ \lambda a + \mu b: a\in A, b\in B \},$$
with $\lambda,\mu\in \R$.

\noindent$\bullet$ The Euclidean distance from a point~$a \in \Rn$
to a set $B\inK$ is defined as
  $$ \mbox{dist}(a,B)=\inf_{b \in B}|a-b|,$$
where $|\cdot|$ is the Euclidean norm in $\Rn$.

\noindent$\bullet$The Hausdorff distance between two sets
$A,B\inK$  is defined by
  $$ \mbox{haus}(A,B)=\max\left\{ \sup_{a \in A}\mbox{dist}(a,B) , \sup_{b
  \in B}\mbox{dist}(b,A) \right\}.$$
\noindent$\bullet$ The set of all projections of $a \in \Rn$ into
a set $B \inK$ is
  $$ \PrjXonY{a}{B}=\{b \in B:|a-b|=\mbox{dist}(a,B)\}.$$

\noindent$\bullet$ For $A,B\inK$ and $0\leq t \leq 1$, the
t-weighted metric average of~$A$ and $B$ is~\cite{Artstein:MA}
  \begin{equation}\label{def_MA}
    \MA{A}{B}{t}=\{ta+(1-t)b:\; (a,b)\in\Pair{A}{B} \}
  \end{equation}
with $\Pair{A}{B} = \{(a,b) \in A \times B:\;
a\in\PrjXonY{b}{A}\;\,
\mbox{or}\;\, b\in\PrjXonY{a}{B} \}$.\\
\noindent The metric average has the metric
property~\cite{Artstein:MA}
\begin{eqnarray}
  \mbox{haus} ( \MA{A}{B}{t}, \MA{A}{B}{s}) &=& |\,t-s|\,\mbox{haus}(A,B), \nonumber\\
  \mbox{haus} ( \MA{A}{B}{t}, A ) &=& (1-t)\,\mbox{haus}(A,B), \label{metric_property}\\
  \mbox{haus} ( \MA{A}{B}{t}, B ) &=& t\,\mbox{haus}(A,B).\nonumber
\end{eqnarray}

\noindent$\bullet$ The modulus of continuity of
$f:[a,b]\rightarrow X$ with images in a metric space $(X,\rho)$ is
\begin{equation}\label{def_ModuliContin}
  \omega_{[a,b]}(f,\delta) = \sup \{\, \rho(f(x),f(y)):\; |x-y|\le \delta,\; x,y \in [a,b] \,\},\;\, \delta >0.
\end{equation}
In this paper $X$ is either $\Rn$ or $\Comp$, and $\rho$ is either
the Euclidean distance or the Hausdorff distance respectively.

\noindent The property of the modulus that we use is
\begin{equation}\label{property_modulus}
  \omega_{[a,b]}(f,\lambda \delta) \le (1+\lambda)\,\omega_{[a,b]}(f,\delta).
\end{equation}

\noindent$\bullet$ By $\Lip{[a,b]}$ we denote the set of all
Lipschitz functions ${f:[a,b]\rightarrow X}$ satisfying
  $$ \rho( f(x), f(y)) \leq {\cal L} |x-y|, \quad \forall\,x,y \in [a,b], $$
where ${\cal L}$ is a constant independent of $x$ and $y.$ \\

\noindent$\bullet$ A variation of ${f:[a,b] \rightarrow X}$ on a
partition ${\chi=\{x_0 <... < x_N:\; x_i \in [a,b],\,i=0,...,N\}}$
is defined by
$$
  V(f,\chi) = \sum_{i=1}^{N} \rho(f(x_i),f(x_{i-1})),
$$
The total variation of $f$ on $[a,b]$ is
$$
  V_{a}^{b}(f) = \sup_{\chi} V(f,\chi).
$$
We say that $f$ is of bounded variation if ${ V_{a}^{b}(f) <
\infty}$ and define in this case
\begin{equation}\label{func_variation}
  v_f(x) = V_{a}^{x}(f),\; x \in [a,b].
\end{equation}
It is obvious that $v_f$ is nondecreasing. If $f$ is also
continuous then $v_f$ is continuous as well. For completeness we
prove it.
\begin{propos}
  A function $f:[a,b] \rightarrow X$ is continuous and of bounded variation
  on $[a,b]$ if and only if $v_f$ is a continuous function on [a,b].
\end{propos}
\begin{proof}
  The sufficiency follows from
  \begin{equation}\label{fContin_if_v(f)Contin}
    \rho(f(x),f(y)) \le V_x^y(f) = v_f(y)-v_f(x),\quad {\rm for}\; x<y.
  \end{equation}
  To prove the other direction, fix $x\in[a,b]$ and $\varepsilon>0$.
  By the uniform continuity of $f$ on $[a,b]$, ${\rho(f(z),f(y))<\varepsilon
  /2}$ if ${|z-y|<\delta}$ for some $\delta>0$. First we show that $v_f$ is
  continuous from the left. We can always choose ${{\chi=\{a=x_0 < x_1 < ... < x_N =
  x\}}}$ such that
  $$
    V_{a}^{x}(f) < V(f,\chi) + \varepsilon/2 = \sum_{i=1}^{N} \rho(f(x_i),f(x_{i-1})) + \varepsilon/2,
  $$
  and ${x-x_{N-1} < \delta}$. Thus
  $$
    V_{a}^{x}(f) < \sum_{i=1}^{N-1} \rho(f(x_i),f(x_{i-1})) + \varepsilon,
  $$
  implying that $v_f(x) - v_f(x_{N-1}) < \varepsilon$.
  By the monotonicity of $v_f$we get for every $x_{N-1}<y<x$
    $$
    v_f(x) - v_f(y) < \varepsilon.
  $$
  Similarly one can show the continuity of $v_f$ from the right.
  Thus we obtain that $v_f$ is continuous at $x$ and consequently
  it is continuous on $[a,b]$.
  \end{proof}

  \noindent From~(\ref{fContin_if_v(f)Contin}) we conclude that
  \begin{equation}\label{moduli_f<moduli Vf}
    {\omega_{[a,b]}(f,\delta)} \le {\omega_{[a,b]}(v_f,\delta)}.
  \end{equation}

\noindent$\bullet$ By CBV we denote the set of all
functions ${f:[a,b]\rightarrow X}$ which are continuous and of bounded variation.\\

\noindent$\bullet$ For a set-valued function $\Map{[a,b]}$, any
single-valued function ${f:[a,b] \rightarrow \Rn}$ with ${f(x) \in
F(x)}$, ${\forall x \in [a,b]}$ is called a selection of $F$.\\

\begin{defin}\label{def_RepresBySelections}
A set of selections of $F$, ${\{f^\alpha:\,\alpha \in {\cal
A}\}}$, is termed a \textbf{representation} of $F$ if
$$
  F(x) = \{ f^\alpha(x):\,\alpha \in {\cal A} \},\quad \forall \,x \in [a,b].
$$
We denote this shortly by ${F=\{f^\alpha:\,\alpha \in {\cal A}
\}.}$
\end{defin}
\bigskip


\section {Linear operators on SVFs based on a \-metric linear combination of ordered sets}
\label{Sect_Chains}

In this section we introduce a new operation on a finite number of
ordered sets. Using this operation we present a new adaptation of
linear operators to multifunctions.

\begin{defin}\label{def_Chain}
  Let ${ \{ A_0,A_1,...,A_N \} }$ be a finite sequence of compact
  sets. A vector ${ (a_0,a_1,...,a_N) }$ with ${a_i \in A_i}$, ${i=0,...,N}$, for which there
  exists $j$, ${0 \leq j \leq N}$ such that
  $$
   a_{i-1} \in \PrjXonY{ a_i }{ A_{i-1} }, \, 1 \leq i \leq j \;\;
   \mbox{and} \;\;
   a_{i+1} \in \PrjXonY{ a_i }{ A_{i+1} }, \, j \leq i \leq {N-1}\,
  $$
  is called a \textbf{metric chain} of $\{A_0, ...,A_N\}$.
\end{defin}
\noindent An illustration of such a metric chain is given in
Figure~\ref{Figure_Chain}.
\begin{center}
    \includegraphics[width=5.43in]{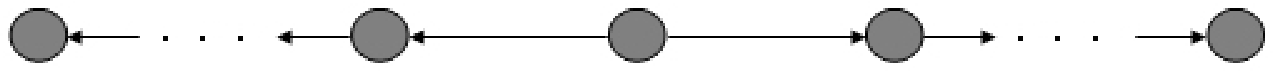}
\begin{footnotesize}
  $$ a_0 \in \PrjXonY{a_1}{A_0} \qquad\quad a_{j-1} \in \PrjXonY{a_j}{A_{j-1}} \qquad\; a_j \in A_j
     \qquad\; a_{j+1} \in \PrjXonY{a_j}{A_{j+1}} \qquad\; a_N \in \PrjXonY{a_{N-1}}{A_N}$$
\end{footnotesize}
\end{center}
\medskip

\begin{pict}\label{Figure_Chain}
  \footnotesize
  {\centering { $\left .\right .$ }\\}
\end{pict}
\medskip

Thus each element of each set $A_i$, ${i=0,...,N}$ generates at
least one metric chain. We denote by ${CH(A_0,...,A_N)}$ the
collection of all metric chains of ${ \{ A_0,...,A_N \} }$. The
set ${CH(A_0,...,A_N)}$ depends on the order of the sets $A_i$,
${i=0,...,N}$.

With this notion of metric chains we can introduce a new operation
between sets.
\begin{defin}
  A \textbf{metric linear combination} of a sequence of sets ${ A_0,...,A_N }$
  with coefficients ${\lambda_0,...,\lambda_N \in \R}$, is
  \begin{equation}\label{def_LinearComb}
    \bigoplus_{i=0}^N \lambda_i A_i =
    \left\{ \sum_{i=0}^N \lambda_i a_i \,: (a_0,...,a_N) \in CH(A_0,...,A_N) \right\}.
  \end{equation}
\end{defin}

Since for two sets $CH(A,B)=\Pair{A}{B}$, in the special case
$N=1$ and ${\lambda_0,\lambda_1 \in [0,1]}$,
${\lambda_0+\lambda_1=1}$, the metric linear combination is the
metric average. The following are two important properties of the
metric linear combination which can be easily seen from the
definition.
\medskip

\noindent $
\begin{array}{ll}
{\rm (i)}& \bigoplus\limits_{i=0}^N\lambda_i A_i =
\bigoplus\limits_{i=0}^N
\lambda_{N-i} A_{N-i}\,, \\
{\rm (ii)} & {\rm For}\;\,  \lambda_0,..., \lambda_N\;\, {\rm
such\; that}\;\, \sum\limits_{i=0}^N \lambda_i = 1,\;\,
\bigoplus\limits_{i=0}^N \lambda_i A = A\,.
\end{array}$
\medskip

With this operation, a large class of linear operators can be
adapted to SVFs.\\
Let $A_{\chi}$, $\chi=\{x_0,...,x_N\}$ be a linear operator of the
form
  \begin{equation}\label{operator_for_f}
    A_{\chi}(f,x) = \sum_{i=0}^N c_i(x)f(x_i),
  \end{equation}
defined on real-valued functions, with domain containing $\chi$.

\begin{defin}\label{def_operator_for_F}
  Let $\Map{[a,b]}$, $\chi \subset [a,b]$ and let
  ${ \{F(x_i), i=0,...,N\} }$ be samples of $F$ at $\chi$. For $A_{\chi}$
  of the form~(\ref{operator_for_f}), we define a \textbf{metric linear operator}
  $A_{\chi}^M$ on $F$ by
  \begin{equation}\label{operator_for_F}
    A^M_{\chi}F(x)=A^M_{\chi}(F,x) = \bigoplus_{i=0}^N c_i(x)F(x_i).
  \end{equation}
We term this operator the \textbf{metric analogue}
of~(\ref{operator_for_f}).
\end{defin}

Note that due to property (ii), the metric analogue of a linear
operator which preserves constants, preserves constant
multifunctions. The analogue of property (ii) does not hold for
Minkowski linear combinations with some negative coefficients,
even for convex sets. This is one reason why only positive
operators, based on Minkowski sum, were applied to set-valued
functions. As is shown in the sequel, Definition
\ref{def_operator_for_F} allows to define also non-positive
operators.

The analysis of the approximation properties of $A_{\chi}^M F$ is
based on properties of the metric piecewise linear approximation
operator.
\bigskip


\section {Metric piecewise linear approximations of SVFs}
\label{Sect_PLApprox}

>From now on $\Map{[a,b]}$, ${\{\, F_i=F(x_i)\,\}_{i=0}^{N}}$,
where ${a=x_0<...<x_N=b}$ and ${\,\chi = (x_0,...,x_N)\,}$ denotes
a partition of $[a,b\,]$. We use the notation
${CH=CH(F_0,...,F_N)}$, and ${\, \delta_{max} = \max\, \{\,
\delta_i:\;\, 0\le i \le N-1 \}}$, ${\; \delta_{min} = \min\, \{\,
\delta_i:\;\, 0\le i \le N-1 \} \,}$ with values ${\delta_i\,}$
defined as ${\,\delta_i = (x_{i+1}-x_i)\,}$, ${\,i=0,...,N-1\,}$.
In case of a uniform partition, we have ${\delta_{max} =
\delta_{min} = h = (b-a)/N}$ and denote such a partition by
$\chi_N$.

\begin{defin}
The \textbf{metric piecewise linear} approximation to $F$ is
$$
    S^M_{\chi}(F,x)= \{\, \lambda_i(x) f_i + (1-\lambda_i(x))f_{i+1}:\, (f_0,...,f_N) \in CH\, \},\quad x\in[x_i,x_{i+1}],
$$
where
  \begin{equation}\label{def_lamda}
    \lambda_i(x) = (x_{i+1}-x)/(x_{i+1}-x_i).
  \end{equation}
\end{defin}

\noindent By construction, the set-valued function $S^M_{\chi}F$
has a representation by selections
  \begin{equation}\label{S_isUnionOf_s}
    S^M_{\chi}F = \{\, s(\chi, \varphi):\, \varphi \in CH(F_0,...,F_N)\,\},
  \end{equation}
where $s(\chi, \varphi)$ is a piecewise linear single-valued
function interpolating the data $(x_i, f_i),$ $i=0,...,N,$ with
$\varphi=(f_0,...,f_N)$.
\bigskip

\noindent Recall the piecewise linear interpolant based on the
metric average, introduced in~\cite{Artstein:MA}:
$$
    S^{MA}_{\chi}(F,x)= \MA{ F_i }{ F_{i+1} }{\lambda_i(x)}, \quad x \in [x_i, x_{i+1}]
$$
with $\lambda_i(x)$ defined by~(\ref{def_lamda}).

It is easy to see by the triangle inequality for the Hausdorff
metric and by~(\ref{metric_property}) that for a continuous
set-valued function $F$
\begin{equation}\label{haus(F,SF)}
  {\rm haus}(F(x),S^{MA}_{\chi}(F,x)) \le 2\,\omega_{[a,b]}(F,\delta_{max}),\; x \in [a,b].
\end{equation}

\begin{remark}\label{remark_S^MA=S^CH}
It is not unexpected that ${S^{MA}_{\chi}F\equiv S^M_{\chi}F}$.

Indeed, for a fixed ${ x \in [x_i, x_{i+1}] }$ and for any ${y \in
S^{MA}_{\chi}(F,x)}$,
$$
  y = \lambda_i(x)f_i + (1-\lambda_i(x))f_{i+1}
$$
with $(f_i,f_{i+1})\in \Pair{F_i}{F_{i+1}}$. Then there exists a
metric chain ${\varphi= (f_0,...,f_i,f_{i+1},...,f_N)}$, ${\varphi
\in CH}$, such that ${y = s(\chi, \varphi)(x)}$. Also it is
obvious that for any ${x \in [a,b]}$ and any ${\varphi \in CH}$,
${s(\chi, \varphi)(x) \in S^{MA}_{\chi}(F,x)}$.
\end{remark}
\medskip

In the following we show that $S^M_{\chi}F$, and its piecewise
linear selections~(\ref{S_isUnionOf_s}) "inherit" some continuity
properties of a continuous multifunction $F$. The following lemma
and corollary consider Lipschitz continuous SVFs.

\begin{lemma}\label{lemma_S_is_Lipshitz}
  Let $F \in \Lip{[a,b]}$, and let $\chi$ be a partition of $[a,b]$. Then
  $$ S^M_{\chi}F \in \Lip{[a,b]} .$$
\end{lemma}
\begin{proof}
For $x,y \in [x_j, x_{j+1}]$ the claim of the lemma follows from
the metric property~(\ref{metric_property}). Now, let ${x \in
[x_j, x_{j+1}]}$ and ${y \in [x_k,x_{k+1}],}$ where ${0 \leq j \le
k \leq N-1.}$ Using the triangle
inequality,~(\ref{metric_property}) and the Lipschitz continuity
of $F$, we get
\begin{align*}
  &\mbox{haus}( S^M_{\chi}(F,x), S^M_{\chi}(F,y))
  \\&\le \frac{x_{j+1}-x}{x_{j+1}-x_j}\,{\rm haus}(F_j, F_{j+1})
  + {\rm haus}(F_{j+1},F_k) + \frac{y-x_k}{x_{k+1}-x_k}\,{\rm haus}(F_k,F_{k+1})
  \\&\le {\cal L} (x_{j+1}-x + x_k-x_{j+1} + y-x_k) \le {\cal L}|y-x|.
\end{align*}
\end{proof}

\begin{corol}\label{corol_selectionsAreLipshitz}
  Under the conditions of Lemma~\ref{lemma_S_is_Lipshitz}, for any
  ${s(\chi,\varphi)}$ in~\ref{S_isUnionOf_s}
  $$ s(\chi,\varphi) \in \Lip{[a,b]} .$$
\end{corol}
\noindent The proof of this corollary is similar to the proof of
the previous lemma and uses the observation that for $k\ge j$
\begin{align*}
  &|s(\chi,\varphi)(x_{j+1}) - s(\chi,\varphi)(x_{k})|
  \le \sum_{l=j+1}^{k-1}|s(\chi,\varphi)(x_l) - s(\chi,\varphi)(x_{l+1})|
  \\&\le \sum_{l=j+1}^{k-1}{\rm haus}(S^M_{\chi}(F,x_l),S^M_{\chi}(F,x_{l+1}))
  \le {\cal L}\sum_{l=j+1}^{k-1}(x_{l+1} - x_l) = {\cal L}|x_k - x_{j+1}|.
\end{align*}
\medskip

\noindent Now we consider the case when $F$ is a general
continuous function.
\begin{lemma}\label{lemma_S_has_ModulusF}
  Let ${\Map{[a,b]}}$ be a continuous set-valued function. Then for any partition
  $\chi$ of $[a,b]$
  $$
     \omega_{[a,b]}(S^M_{\chi}F,\delta) \le 5\,\omega_{[a,b]}(F,\delta).
  $$
\end{lemma}
\begin{proof}
By definition, for any $\delta > 0$
$$
  \omega_{[a,b]}(S^M_{\chi}F,\delta) = \sup \,\{\, {\rm haus}(S^M_{\chi}(F,x),S^M_{\chi}(F,y)):\;
  |x-y|\le \delta,\; x,y \in [a,b] \,\}.
$$

In case ${x,y \in [x_j,x_{j+1}]}$, ${|x-y|\le \delta}$, the claim
of the lemma is obtained using~(\ref{def_lamda}), the metric
property~(\ref{metric_property}) and~(\ref{property_modulus}),
\begin{equation}\label{modulus_resultA}
    \begin{split}
 {\rm haus}(S^M_{\chi}(F,x),S^M_{\chi}(F,y)) &= \frac{|x-y|}{\delta_j}\,{\rm haus}(F_j,F_{j+1})
  \le \frac{\delta}{\delta_j}\omega_{[a,b]}(F,\delta_j)
  \\&\le \frac{\delta}{\delta_j}\left(1 + \frac{\delta_j}{\delta}\right )\omega_{[a,b]}(F,\delta)
  \le 2\,\omega_{[a,b]}(F,\delta).
  \end{split}
  \end{equation}

Now, let ${x \in [x_j, x_{j+1}]}$, ${y \in [x_k,x_{k+1}],}$ ${0
\leq j < k \leq N-1}$ and ${|x-y|\le \delta}$. By the triangle
inequality
\begin{equation}\label{modulus_resultB}
\begin{split}
  {\rm haus}(S^M_{\chi}(F,x),S^M_{\chi}(F,y)) &\le {\rm haus}(S^M_{\chi}(F,x),S^M_{\chi}(F,x_{j+1}))
  \\&+ {\rm haus}(S^M_{\chi}(F,x_{j+1}),S^M_{\chi}(F,x_k))
  \\&+ {\rm haus}(S^M_{\chi}(F,x_k),S^M_{\chi}(F,y)),
\end{split}
\end{equation}
while by the interpolation property of $S^M_{\chi}F$ and since
$|x_k-x_{j+1}|\le \delta$, we have
\begin{equation}\label{modulus_resultC}
  {\rm haus}(S^M_{\chi}(F,x_{j+1}),S^M_{\chi}(F,x_k)) \le \omega_{[a,b]}(F,\delta).
\end{equation}
Applying~(\ref{modulus_resultA}) and~(\ref{modulus_resultC})
to~(\ref{modulus_resultB}) we obtain the claim of the lemma.
\end{proof}
\smallskip

\begin{corol}\label{corol_selections_LocalModulusS}
  For any $s(\chi,\varphi)$ in~(\ref{S_isUnionOf_s}) and any $x,y \in
  [x_j,x_{j+1}]$, $0\le j\le N-1$
  \begin{equation}\label{dist_s<moduli_F}
     |s(\chi,\varphi)(x)-s(\chi,\varphi)(y)| \le 2\,\omega_{[a,b]}(F,|x-y|)
  \end{equation}
  Also,
  \begin{equation}\label{moduli_s<moduli_F}
     \omega_{[a,b]}(s(\chi,\varphi),\delta) \le 4\,\omega_{[a,b]}(F,\delta),\quad \delta \le
     \delta_{min}\,.
  \end{equation}
\end{corol}
\noindent   The proof of this corollary is similar to the proof of
assertion~(\ref{modulus_resultA}).
\bigskip

We cannot generalize~(\ref{moduli_s<moduli_F}) for arbitrary
${0<\delta \le b-a}$ if $F$ is only continuous. Yet we can get an
estimate for $\omega_{[a,b]}(s(\chi,\varphi),\delta)$ if $F$ is
continuous and of bounded variation.

\begin{lemma}\label{lemma_omega(s)<omega(v(F))}
Let $F \in CBV([a,b])$. Then for any $s(\chi,\varphi)$
in~(\ref{S_isUnionOf_s}),
  $$
     \omega_{[a,b]}(s(\chi,\varphi),\delta)
     \le 4\,\omega_{[a,b]}(F,\delta) + \omega_{[a,b]}(v_F,\delta)
    \le 5\,\omega_{[a,b]}(v_F,\delta).
  $$
\end{lemma}
\begin{proof}
  Denote $s=s(\chi,\varphi)$. For a given $\delta >0$, let ${x \in [x_j, x_{j+1}]}$, ${y \in [x_k,x_{k+1}]}$,
  ${0\leq j \le k \leq N-1}$, such that ${|x-y|\le \delta}$. Using the definition of $s(\chi,\varphi)$
  and of $S^M_{\chi}F$ we get
  \begin{align*}
    &|s(x)-s(y)| \le |s(x)-s(x_{j+1})| + \sum_{l=j+1}^{k-1} |s(x_{l+1})-s(x_l)| + |s(y) -s(x_k)|
    \\&\le \frac{x_{j+1}-x}{\delta_j}|s(x_{j+1})-s(x_j)| + \sum_{l=j+1}^{k-1} {\rm haus}(F(x_{l+1}), F(x_l))
    + \frac{y-x_k}{\delta_k}|s(x_{k+1})-s(x_k)|\,.
  \end{align*}
  Now,~(\ref{dist_s<moduli_F}) yields
  $$
    |s(x)-s(y)| \le 4\,\omega_{[a,b]}(F,\delta) + V_{x_{j+1}}^{x_k}(F) \le 4\,\omega_{[a,b]}(F,\delta) +
    \omega_{[a,b]}(v_F,\delta).
  $$
Taking the supremum over ${|x-y|\le \delta}$ and
using~(\ref{moduli_f<moduli Vf}), we complete the proof.
\end{proof}
\bigskip


\section {Approximation by metric linear operators}
\label{Sect_ApproxByOperators}

We use the metric piecewise linear approximation to obtain error
estimates for metric linear operators.

Let $A^M_{\chi}F$ be defined by~(\ref{operator_for_F}) and
$S^M_{\chi}F$ be a metric piecewise linear multifunction as
defined in Section~\ref{Sect_PLApprox}. By
Definition~\ref{def_operator_for_F}
\begin{equation}\label{AF=A_SF}
  A^M_{\chi}F \equiv A^M_{\chi}(S^M_{\chi}F).
\end{equation}
Moreover by~(\ref{operator_for_f}),~(\ref{operator_for_F})
and~(\ref{S_isUnionOf_s})
\begin{equation}\label{AS=UnionOf_As}
    A^M_{\chi}(S^M_{\chi}F) = \left\{ A_{\chi}s(\chi,\varphi):\, \varphi \in CH(F_0,...,F_N)\right\}.
\end{equation}

\begin{remark}\label{remark_AF=UnionOf_As}
In contrast to our previous definition of positive operators for
SVFs based on the metric average~\cite{Dyn-Farkhi:Subdiv_MA,
Dyn-Mokhov:Approx_MA}, the metric analogues~(\ref{operator_for_F})
of two linear operators of the form~(\ref{operator_for_f}), which
are identical on single-valued functions, are identical on SVFs.
For example, in~\cite{Dyn-Farkhi:Subdiv_MA, Dyn-Mokhov:Approx_MA}
spline subdivision schemes are not identical to the Schoenberg
spline operators for SVFs.
\end{remark}
\medskip

The metric analogues of linear operators of the
form~(\ref{operator_for_f}), which approximate real-valued
functions, are approximating SVFs.
By~(\ref{AF=A_SF}),~(\ref{AS=UnionOf_As}) the approximation
results depend on the way $A_{\chi}$ approximates piecewise linear
real-valued functions.

In what follows ${\phi:[a,b]\times \R_+ \rightarrow \R_+}$ is a
continuous real-valued function, non-decreasing in its second
argument, satisfying ${\phi(x,0)= 0}$, and ${\cal S}_{\chi}$
denotes the set of piecewise linear continuous single-valued
functions, with values in $\Rn$ and knots at $\chi$.

\begin{theorem}\label{ApproxTheorem_FisLipshitz}
  Let $A_{\chi}$ be of the form~(\ref{operator_for_f}), such that for any ${s \in {\cal S}_{\chi}\bigcap\Lip{[a,b]}}$
  \begin{equation}\label{dist(As,s)_Lipschitz}
    |A_{\chi}(s,x)-s(x)| \leq C\,{\cal L}\phi(x,\delta_{max}).
  \end{equation}
  \noindent Then if $F \in \Lip{[a,b]}$,
  $$
    {\rm haus}(A^M_{\chi}(F,x),F(x)) = 2\,{\cal L}\delta_{max} + C {\cal L}\phi(x,\delta_{max}).
  $$
\end{theorem}
\begin{proof}
  By~(\ref{AF=A_SF})
  \begin{equation}\label{theoremLip_triangle}
  \begin{split}
     {\rm haus}(A^M_{\chi}(F,x),F(x)) \le {\rm haus}(A^M_{\chi}(S^M_{\chi}F,x),S^M_{\chi}(F,x))
     + {\rm haus}(S^M_{\chi}(F,x),F(x)),
  \end{split}
  \end{equation}
  while by~(\ref{AS=UnionOf_As})
  $$
      {\rm haus}( A^M_{\chi}(S^M_{\chi}F,x), S^M_{\chi}(F,x) ) \leq
      \sup_{\varphi \in \, CH}
      |A_{\chi}(s(\chi,\varphi),x)-s(\chi,\varphi)(x)|\,.
  $$
  \noindent In view of
  Corollary~\ref{corol_selectionsAreLipshitz} and~(\ref{dist(As,s)_Lipschitz})
  \begin{equation}\label{sup_dist(As,S)_bounded}
    \sup\limits_{\varphi \in CH} |A_{\chi}(s(\chi,\varphi),x)-s(\chi,\varphi)(x)| \le C \,{\cal L}\phi(x,\delta_{max})
  \end{equation}
  \noindent The proof is completed by
  substituting~(\ref{sup_dist(As,S)_bounded})
  and~(\ref{haus(F,SF)}) in~(\ref{theoremLip_triangle}).
\end{proof}
\medskip

For general continuous SVFs we cannot prove an analogous
approximation result. Yet for continuous multifunctions of bounded
variation we get a weaker approximation result, by applying
Lemma~\ref{lemma_omega(s)<omega(v(F))} instead of
Corollary~\ref{corol_selectionsAreLipshitz} in the proof of
Theorem~\ref{ApproxTheorem_FisLipshitz}.

\begin{theorem}\label{ApproxTheorem_FisCont+BV}
Let $F \in CBV([a,b])$, and let $A_{\chi}$ be of the
form~(\ref{operator_for_f}), satisfying
  \begin{equation}\label{dist(As,s)_Cont+BV}
    |A_{\chi}(s,x)-s(x)| \leq
    C\,\omega_{[a,b]}(s,\phi(x,\delta_{max})),\quad s\in {\cal S}_{\chi}\,.
  \end{equation}
Then
  $$
    {\rm haus}(A^M_{\chi}(F,x),F(x)) = 2\,\omega_{[a,b]}(F,\delta_{max}) + 5C\omega_{[a,b]}(v_F,\phi(x,\delta_{max})).
  $$
\end{theorem}
\medskip

For continuous SVFs which are not of bounded variation we can
prove an approximation result only for uniform partitions and for
a limited class of linear operators.

\begin{theorem}\label{ApproxTheorem_FisContinuous}
  Let $A_N$ be a linear operator of the
  form~(\ref{operator_for_f}), defined on a uniform partition $\chi_N$ and let
  ${h=(b-a)/N}$. If
  \begin{equation}\label{dist(As,s)_Continuous}
    |A_N(s,x)-s(x)| \leq C\,\phi(x,\omega_{[a,b]}(s,h)),\quad s \in {\cal S}_{\chi}\,,
  \end{equation}
  then for a continuous $F$
  $$
    {\rm haus}(A^M_N(F,x),F(x)) = 2\,\omega_{[a,b]}(F,h) + C \phi(x,4\,\omega_{[a,b]}(F,h)).
  $$
\end{theorem}
The proof of this result repeats the proof of
Theorem~\ref{ApproxTheorem_FisLipshitz}, but replaces
Corollary~\ref{corol_selectionsAreLipshitz}
by~(\ref{moduli_s<moduli_F}) of
Corollary~\ref{corol_selections_LocalModulusS}.
\bigskip


\section {Examples}\label{Sect_Examples}

In this section we present metric analogues for SVFs of the
Schoenberg spline operators and the Bernstein polynomial operators
and give approximation results. We conclude by two examples
demonstrating the operation of metric analogues of parabolic
interpolants. To our knowledge so far only positive operators were
applied to SVFs. The two examples we present assert that such
interpolation between sets is reasonable.

\subsection {Metric Bernstein operators}\label{Bernstein}

The Benstein operator $\Bern{f}$ for $f \in C[0,1]$ is
\begin{equation}\label{BernScalar}
  \Bern{f} = \sum_{i=0}^{N} {N \choose i} x^i (1-x)^{N-i} f \left( \frac{i}{N} \right ).
\end{equation}
It is known (see \cite{DeVore_Lorentz:93}, Chapter 10) that there
exists a constant $C$ independent of $f$ such that
\begin{equation}\label{dist_BernScalar}
  |f(x)-\Bern{f}| \le C \omega_{[0,1]}(f,\sqrt{x(1-x)/N}).
\end{equation}
The classical Bernstein operator for $\Map{[0,1]}$ with sums of
numbers replaced by Minkowski sums of sets is
\begin{equation}\label{BernMink}
  B_N^{Mn}(F,x)= \sum_{i=0}^{N} {N \choose i} x^i (1-x)^{N-i} F \left(
  \frac{i}{N} \right ).
\end{equation}
It was shown in~\cite{Vitale:Bern_Mink} that for ${x \in(0,1)}$
the limit of ${B_N^{Mn}(F,x)}$ when ${N \to \infty,}$ is the
convex hull of~$F(x)$, therefore these operators cannot
approximate SVFs with general images.

In~\cite{Dyn-Mokhov:Approx_MA} Bernstein operators for set-valued
functions are defined procedurally in terms of the de~Casteljau
algorithm, with the metric average as a basic binary operation,
\begin{eqnarray}
  F_i^0 = F(i/N),\quad i=0,...,N,\qquad\qquad\qquad\;\nonumber\\
  F_i^k = \MA {F_i^{k-1}}{F_{i+1}^{k-1}}{1-x},\quad  i=0,1,...,N-k, \; k=1,...,N, \label{Bern_MA_Recur}\\
  B_N^{MA}(F,x) = F_0^N.\qquad\qquad\qquad\qquad\nonumber
\end{eqnarray}

We do not know whether these operators approximate multifunctions
with general compact images in $\Rn$, yet they approximate
multifunctions with compact images in $\R$ all consisting of the
same number of disjoint intervals~\cite{Dyn-Mokhov:Approx_MA}.
\bigskip

\noindent Here we investigate the metric analogue of the Bernstein
operators for SVFs.
\begin{defin}\label{def_Bern_Chains}
  For $\Map{[0,1]}$ the metric Bernstein operator is
  \begin{align*}
    B^M_N(F,x) &= \bigoplus_{i=0}^N {N \choose i} x^i (1-x)^{N-i}F\left(\frac{i}{N}\right)
    \\&= \left\{ \sum_{i=0}^N {N \choose i} x^i (1-x)^{N-i} f_i\,:\; (f_0,...,f_N) \in CH
    \right\},
  \end{align*}
  where $CH=CH(F(0),F(1/N)...,F(1))$.
\end{defin}

\noindent By Theorem~\ref{ApproxTheorem_FisLipshitz} and
by~(\ref{dist_BernScalar}) we conclude that
\begin{corol}
  Let $F\in\Lip{[0,1]}$, then
 $${\rm haus}(B_N^M(F,x),F(x)) \le 2{\cal L}/N + C{\cal L}\sqrt{x(1-x)/N}.$$
\end{corol}

\noindent Moreover by Theorem~\ref{ApproxTheorem_FisCont+BV} and
by~(\ref{dist_BernScalar})
\begin{corol}
  Let $F \in CBV([0,1])$, then
 $${\rm haus}(B_N^M(F,x),F(x)) \le 2\,\omega_{[0,1]}(F,1/N) + 5C\omega_{[0,1]}(v_F,\sqrt{x(1-x)/N}).$$
\end{corol}

\noindent Since~(\ref{dist(As,s)_Continuous}) does not hold for
these operators, Theorem~\ref{ApproxTheorem_FisContinuous} cannot
be applied.
\bigskip

\subsection {Metric Schoenberg operators}\label{Schoenberg}

For a uniform partition $\chi_N$, the "classical" set-valued
analogues of the Schoenberg spline operators for $\Map{[0,1]}$ is
\begin{equation}\label{SchoenScalar}
  S^{Mn}_{m,N}(F,x)=\Sum{i=0}{N}F(i/N)\Bspline{m}{Nx-i},
\end{equation}
where $\Bspline{m}{x}$ is the B-spline of order~$m$ (degree $m-1$)
with integer knots and support $[0,m]$, and where the linear
combination is in the Minkowski sense. An example, given
in~\cite{Vitale:Bern_Mink}, shows that the operators
in~(\ref{SchoenScalar}) with $m=2$ and $N\rightarrow \infty$
cannot approximate $F$ with general compact images, in any point
of ${[0,1]\setminus \chi_N}$.

A Shoenberg operator based on the metric average is introduced
in~\cite{Dyn-Mokhov:Approx_MA}, by a procedural definition in
terms of repeated binary averages according to the de~Boor
algorithm. It is proved that for H\"older continuous set-valued
functions, the approximation rate is the H\"older exponent.

Here we consider the metric analogue of the Schoenberg operators.
\begin{defin}\label{def_Spline_Chains}
  The metric Shoenberg operator of order $m$ for a set-valued function ${\Map{[0,1]}}$
  and a uniform partition ${\chi_N}$ is defined by
  \begin{align*}
    S_{m,N}^M(F,x) = \bigoplus_{i=0}^N \Bspline{m}{Nx-i} F\left(\frac{i}{N}\right)
    = \left\{ \sum_{i=0}^N \Bspline{m}{Nx-i} f_i:\; (f_0,...,f_N) \in CH
    \right\},
  \end{align*}
  where $CH=CH(F(0),F(1/N)...,F(1))$.
\end{defin}

By Theorem~\ref{ApproxTheorem_FisContinuous} and the known
approximation result in case of single-valued functions~(see
\cite{deBoor:Splines}, Chapter XII), we obtain
\begin{corol}
  Let $F$ be a continuous SFV defined on $[0,1]$. Then
$$
  {\rm haus}(S_{m,N}^M(F,x), F(x)) = 2 \left( 1+2 \left\lfloor \frac{m+1}{2}\right\rfloor \right)
  \omega_{[0,1]}(F,1/N), \quad x \in \left [\frac{m-1}{N},1\right ]
$$
with $\lfloor t \rfloor$ the maximal integer not greater than $t$.
\end{corol}

The approximation result in the specific case of Lipschitz
continuous SVFs, can be further improved by applying
Theorem~\ref{ApproxTheorem_FisLipshitz}.
\begin{corol}
  For $F \in \Lip{[0,1]}$,
$$
  {\rm haus}(S_{m,N}^M(F,x), F(x)) = \left(2+ \left\lfloor \frac{m+1}{2}\right\rfloor \right)
  \frac{{\cal L}}{N}.
$$
\end{corol}
\bigskip

\subsection {Metric polynomial interpolants}\label{Interpolation}

\begin{defin}\label{def_Bern_Chains}
  Let $\Map{[a,b]}$, and let $\chi$ be a partition of $[a,b]$.
  The metric polynomial interpolation operator is given by
  \begin{align*}
    P_{\chi}^M(F,x) = \bigoplus_{i=0}^N l_i(x) F(x_i)
    = \left\{ \sum_{i=0}^N l_i(x) f_i\,:\, (f_0,...,f_N) \in CH(F(x_0),...,F(x_N))
    \right\},
  \end{align*}
  with $l_i(x)$ the $i$-th Lagrange polynomial,
  $$l_i(x) = \prod\limits_{j=0,j\neq i}^N \frac{x-x_j}{x_i-x_j}\,.$$
\end{defin}
\smallskip

 To illustrate our method we apply the metric parabolic interpolation
 operator to three sets in~$R$. We consider two different
 examples.

  \noindent The first example: $x_0 = 0$, $x_1 = 2$, $x_2 = 6$;
  $$
  F(x_0) = [2,8],\qquad  F(x_1) = \{5\},\qquad F(x_2) = \{5\}.
  $$
  The second example: $x_0 = 0$, $x_1 = 4$, $x_2 = 8$;
  $$
  F(x_0) = [2,4]\cup[6,8],\qquad  F(x_1) = [4.5,5.5],\qquad F(x_2) =
  [2,4]\cup[6,8].
  $$
The two set-valued interpolants are illustrated in
Figure~\ref{Figure_Fish} and Figure~\ref{Figure_Butterfly}
respectively.

\begin{center}
    \includegraphics[width=5.25in]{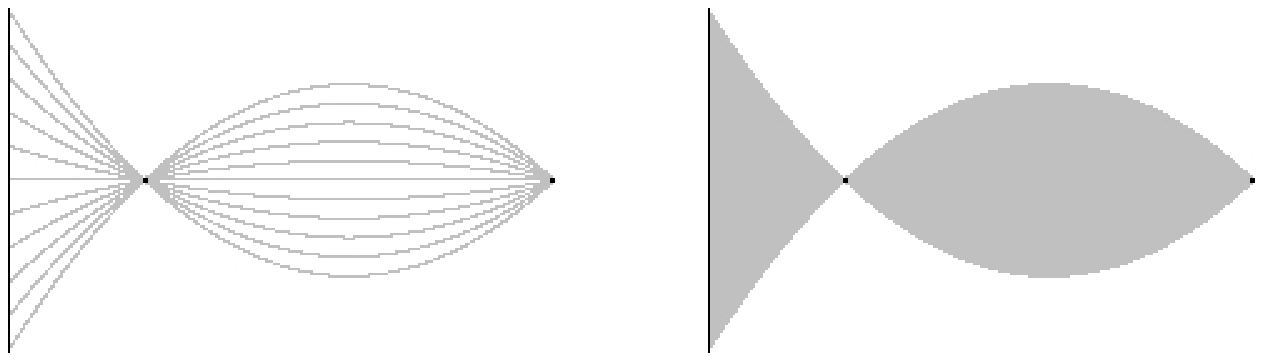}
\begin{footnotesize}
\end{footnotesize}
\end{center}
\begin{pict}\label{Figure_Fish}
  {\centering { Parabolic interpolation - first example. }\\}
\end{pict}

\begin{center}
    \includegraphics[width=5.4in]{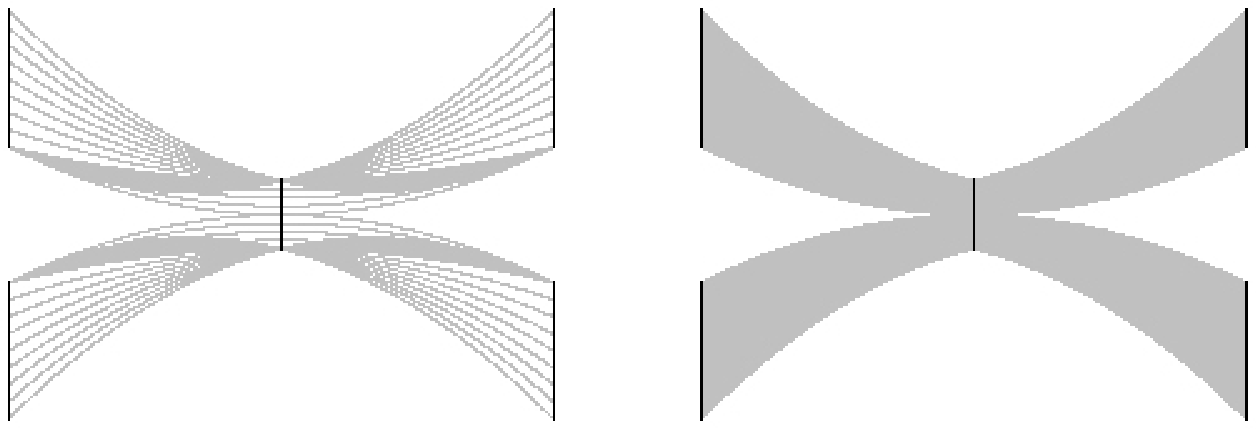}
\begin{footnotesize}
\end{footnotesize}
\end{center}
\begin{pict}\label{Figure_Butterfly}
  {\centering { Parabolic interpolation - second example. }\\}
\end{pict}
In the above figures the sets in black are $F(x_0)$, $F(x_1)$,
$F(x_2)$ and the gray curves are the parabolic interpolants to the
selections in~(\ref{S_isUnionOf_s}).


\begin{thebibliography}{99}

\bibitem{Artstein:MA}
  Z. Artstein,
  Piecewise linear approximations of set-valued maps,
  Journal of Appro\-ximation Theory 56, (1989), pp.~$41$-$47$.

\bibitem{Dommisch:LipSelections}
  G. Dommisch,
  On the existence of Lipschitz continuous and differentiable
  selections for multifunctions, in Parametric Optimization and
  Related Topics, J.Guddat, H.Th.Jongen, B.Kummer, F.No\v{z}i\v{c}ka (eds.),
  Akademie-Verlag Berlin,~(1987), pp.~$60$-$73$.

\bibitem{deBoor:Splines}
  C. de Boor,
  A practical guide to spline. Springer-Verlag, (2001).

\bibitem{DeVore_Lorentz:93}
  R.~DeVore and G.~Lorentz, Constructive Approximation,
  Springer-Verlag Berlin, (1993).

\bibitem{Dyn-Farkhi:Subdiv_MA}
  N.Dyn, E.Farkhi,
  Spline subdivision schemes for compact sets with metric
  averages, in Trends in Approximation Theory, K.Kopotun, T.Lyche
  and M.Neamtu (eds.), Vanderbilt Univ. Press, (2001), pp.~$95$-$104$ .

\bibitem{Dyn-Farkhi:ConvRates}
  N.Dyn, E.Farkhi,
  Set-valued approximations with Minkowski
  averages - convergence and convexification rates, Numerical
  Functional Analysis and Optimization 25, (2004), pp.~$363$-$377$.

\bibitem{Dyn-Mokhov:Approx_MA}
  N.~Dyn and A.~Mokhov,
  Approximations of set-valued functions based
  on the metric average, submitted.

\bibitem{HenryHermes:ContSelections}
  H. Hermes, On continuous and measurable selections and the
  existence of solutions of genera\-lized differential equations,
  Proceeding of the American Mathematical Society 29, (1971),
  pp.~$535$-$542$.

\bibitem{Vitale:Bern_Mink}
  R.A. Vitale,
  Approximation of convex set-valued functions,
  Journal of Approximation Theory 26, (1979), pp.~$301$-$316$.

\bibitem{DynWal}
  J. Wallner and N. Dyn, Convergence and $C^1$ analysis of
  subdividision schemes on manifolds by proximity, to appear in
  CAGD.

\end{thebibliography}
\end{document}